\documentclass[11pt]{article}
\usepackage{amsthm}
\usepackage{amsmath}
\usepackage{amssymb}
\usepackage{graphicx}
\usepackage{enumerate}
\usepackage[pdfstartview=FitH]{hyperref}
\hypersetup{
    pdftitle=   {Large complete minors in expanding graphs } pdfauthor ={ Younjin Kim}}

\newtheorem{definition}{Definition}[section]

\newtheorem{theorem}[definition]{Theorem}

\newtheorem{lemma}[definition]{Lemma}

\numberwithin{equation}{section}

\providecommand{\keywords}[1]
{
  \small	
  \textbf{{Keywords}} #1
}

\newcommand{\ignore}[1]{}

\if11     % 11 = on,  10 = off
\usepackage[mathlines]{lineno}
\newcommand*\patchAmsMathEnvironmentForLineno[1]{%
  \expandafter\let\csname old#1\expandafter\endcsname\csname #1\endcsname
  \expandafter\let\csname oldend#1\expandafter\endcsname\csname end#1\endcsname
  \renewenvironment{#1}%
     {\linenomath\csname old#1\endcsname}%
     {\csname oldend#1\endcsname\endlinenomath}}% 
\newcommand*\patchBothAmsMathEnvironmentsForLineno[1]{%
  \patchAmsMathEnvironmentForLineno{#1}%
  \patchAmsMathEnvironmentForLineno{#1*}}%
\AtBeginDocument{%
\patchBothAmsMathEnvironmentsForLineno{equation}%
\patchBothAmsMathEnvironmentsForLineno{align}%
\patchBothAmsMathEnvironmentsForLineno{flalign}%
\patchBothAmsMathEnvironmentsForLineno{alignat}%
\patchBothAmsMathEnvironmentsForLineno{gather}%
\patchBothAmsMathEnvironmentsForLineno{multline}%
}
\fi

\def\inst#1{$^{#1}$}

\begin{document}

\title{Large complete minors in expanding graphs}

\author{ Younjin Kim\inst{1}\thanks{The author was supported by Basic Science Research Program through the National Research Foundation of Korea (NRF) funded by the Ministry of Education (2017R1A6A3A04005963).}}

\maketitle

\begin{center}
{\footnotesize
\inst{1} 
Institute of Mathematical Sciences, Ewha Womans University, Seoul, South Korea \\
\texttt{younjinkim@ewha.ac.k}r
\\\ \\
}
\end{center}

\begin{abstract}
 
In 2009, Krivelevich and Sudakov  studied the existence of large  complete minors in $(t,\alpha)$-expanding graphs whenever the expansion factor $t$ becomes super-constant. In this paper, we  give an extension of the  results of Krivelevich and Sudakov
by investigating  a connection between the existence of large complete minors in graphs and good vertex expansion properties.  \\


\end{abstract}

\section{Introduction}

A graph here will refer to a simple undirected graph without loops. An undirected graph $G$ is a set $V=V(G)$ of vertices and a set $E(G)$ of edges of $G$. A graph $H$ is a {\it minor} of
$G$, denoted by $G \succ H$, if it can be obtained from $G$ by vertex deletions, edge deletions, and contractions.
A graph $G$ is an {\it $\alpha$-expander} if every set $A$ of at most half of the vertices of $G$ has at least $\alpha\cdot|A|$
outside neighbors in $G$.  A graph $G$ is a {\it $(t,\alpha)$-expanding} if every subset $X \subset V(G)$ of size $|X| \leq \frac{\alpha |V(G)|}{t}$ has at least $t|X|$ external neighbors in $G$. 
An abundance of various results and conjectures providing sufficient conditions for the existence of large complete minors have been widely investigated in the context of undirected simple graphs.  For example, the famous Hadwiger's conjecture addressed a connection between the chromatic number and the existence of large minors, which states if $G$ is a graph with the chromatic number $\chi(G)$, then $G$ contains a clique minor of size at least $\chi(G)$. Other authors also discussed a connection between the existence of large minors and graphs with large girth~\cite{Diestel2004, Krivelevich2009,Kuhn2003, Thomassen1983}, $K_{s,t}$-free graphs~\cite{Krivelevich2009,Kuhn2004}, graphs without small vertex separators~\cite{Alon1990,Kawar2010,Plotkin1994}, lifts of graphs~\cite{Drier2006}, random graphs~\cite{Fount2008}, random regular graphs\cite{Fount2009}, and others.\\

  The existence of large  complete minors in $(t,\alpha)$-expanding graphs with the expansion factor $t=\Theta(1)$ was studied by 
Alon, Seymour, and Thomas~\cite{Alon1990}, Plotkin, Rao, and Smith~\cite{Plotkin1994}, and  Kleinberg and Rubinfeld~\cite{Kleinberg1996}. Especially, in 
  $1996$, Kleinberg and Rubinfeld~\cite{Kleinberg1996}  proved that for every fixed $\alpha > 0$, there exists a constant $c > 0$ such that an  $\alpha$-expander graph of order $n$ contains every graph $H$ with
at most $\frac{n}{\log^cn}$ vertices and edges as a minor. As an extension of their result, in 2009, Krivelevich and Sudakov~\cite{Krivelevich2009} obtained a connection between  complete minors in  graphs and vertex expansion properties as follows.\\

\begin{theorem}[Krivelevich, Sudakov~\cite{Krivelevich2009}]
Let $0 < \alpha < 1 $ be a constant, and let $t \geq 10$.   Then a $(t,\alpha)$-expanding graph of order $n$ contains a minor with an average degree of at least

$$ c\frac{\sqrt{nt\log t}}{\sqrt{\log n}},$$

\noindent where $c=c(\alpha) > 0$ is a constant.\\
\end{theorem}

 With the results of Kostochka~\cite{Kostochka1982} and Thomason~\cite{Thomason1984}, Krivelevich and Sudakov~\cite{Krivelevich2009} also got the following result, which  enables us to show the existence of larger minors in $(t,\alpha)$-expanding graphs whenever the expansion factor $t$ becomes super-constant.

\begin{theorem}[Krivelevich, Sudakov~\cite{Krivelevich2009}]
 Let $0 < \alpha < 1 $ be a constant, and let $t \geq 10$.   Then  a $(t, \alpha)$-expanding graph of order $n$ contains a clique minor of size 

 $$ \Omega \left(\alpha^2 \frac{\sqrt{nt\log t}}{\sqrt{\log n}}\right).$$\\
\end{theorem}

  Recently Krivelevich and Nenadov~\cite{Krivelevich2019} complemented the results of Krivelevich and Sudakov~\cite{Krivelevich2009} by investigating a connection  between the contraction clique number of a graph and  good edge expansion properties.  Usually, an advantage of edge expansion  over vertex expansion is that it is easier to verify. In this paper, we  investigate a connection between complete minors in graphs and good vertex expansion properties as follows.
In the following theorem, we use the expander notation introduced by 
  Haslegrave, Kim, and Liu~\cite{HKL20} such that similar expansion occurs even after removing a relatively small set of edges.
 Our result is   an extension of the  results of Krivelevich and Sudakov~\cite{Krivelevich2009}. An  $n$-vertex graph $G$ with the minimum degree $d$ is {\it $\epsilon$-locally sparse} if for every subset $U \subseteq V(G)$ of size $|U| \leq \epsilon n$, then an average degree $d(G[U])$ on the induced subgraph $G[U]$ is  at most $ \epsilon d$. An $n$-vertex graphs $G$ is  {\it $\epsilon$-vertex-expand} if for every subset $U \subseteq V(G)$ of size $|U|\leq \frac{n}{2}$, we have $|N(U)|\geq \epsilon|U|$.\\
 
  \begin{theorem}\label{maintheorem} 
  For $\epsilon > 0$ and $d > 0$, we let $G$ be a graph of order $n$ satisfying the following conditions:
  \begin{itemize}
  \item $\frac{\epsilon}{2}$-locally sparse
  \item $\epsilon$-vertex-expand
  \item minimum degree $d$.
  \end{itemize}
  Then  $G$ contains a clique minor of size 

 $$ \tilde{\Omega}\left(\frac{\sqrt{nd}}{\log^{10}n}\right).$$\\

\end{theorem}

Our proof utilizes robust sublinear expanders, which is an extension of a notion of expander introduced in the $90$s by Koml\' os and Szemer\' edi~\cite{KS96}. It has proven to be a powerful tool for embedding sparse graphs, playing an essential role in the recent resolution of several long-standing conjectures that were previously out of reach.(see \cite{FKKL2021}, \cite{LM2017}, \cite{LM2021})\\

Our paper is organized as follows. In Section $2$, 
 we introduce some preliminaries on expanders. 
In Section $3$, we give the proof of  Theorem~\ref{maintheorem}.

 \section{Preliminaries}
 We denote by $[n]$ the set $\{1,2,\cdots, n\}$ of the first $n$ positive integers. Given a set $X$ and $k\in \mathbb{N}$, we let  ${ {X}\choose {k}}$ for the family of its $k$-element subsets. For brevity, we write $v$ for a singleton set $\{ v \}$ and $xy$ for a set of pairs $\{x, y \}$. 
Given a graph $G$, denote its average degree $2e(G)/|G|$ by $d(G)$. Let $F\subseteq G$ and $H$ be graphs, and $U\subseteq V(G)$. We write $G[U]\subseteq G$ for the induced subgraph of $G$ on vertex set $U$. Denote by $G\cup H$  the graph with vertex set $V(G)\cup V(H)$ and edge set $E(G)\cup E(H)$, and write $G-U$ for the induced subgraph $G[V(G)\setminus U]$, and $G\setminus F$ for the spanning subgraph of $G$ obtained from removing the edge set of $F$. 

\subsection{Sublinear expander}
Expanders are typically well-connected sparse graphs in which vertex subsets exhibit expansions. Originally introduced for network design, expanders, apart from being a central notion in graph theory, have also close interplay with other areas of mathematics and as well as theoretical computer science. This is partially reflected by the fact that expanders have equivalent definitions from different angles. Indeed, in terms of graph expansions, an expander is a graph whose vertex subsets have a large neighborhood: while algebraically, there is a spectral gap when you look at its adjacency matrix; and probabilistically, random walks on expanders are rapidly mixing.\\

In 1994,  Koml\'os and Szemeredi~\cite{KS96} first introduced the following notion of graph expansion. For $\epsilon_1 >  0$ and $ k > 0$, we let $\rho(x)$ be the function

$$\rho(x) = \rho(x,\epsilon_1,t) : =  \left\{ \begin{array}{lcl}
 0 & \mbox{if} & x < t/5, \\
\epsilon_1/\log^2(15x/t) & \mbox{if} & x \geq t/5.\\
\end{array}\right.$$\\

\noindent  Note that when $x \geq t/2$,  $\rho(x)$ is decreasing, while $\rho(x)\cdot x$ is increasing. 
Koml\'os and Szemeredi~\cite{KS96} said that a graph $G$ is an ($\epsilon_1,t$)-expander if for any subset $X \subseteq V(G)$ of size $t/2 \leq |X| \leq |G|/2$, we have $ |N_G(X)| \geq \rho(|X|)|X|.$\\

\noindent Recently, Haslegrave, Kim, and Liu~\cite{HKL20} extended the expander notation  such that similar expansion occurs even after removing a relatively small set of edges as follows. We will use the following sublinear expander introduced by Haslegrave, Kim, and Liu~\cite{HKL20}.

\begin{definition} [($\epsilon_1,t$)-expander] A graph $G$ is an ($\epsilon_1,t$)-expander if 
\begin{itemize}
\item for any subset $X \subseteq V(G)$ of size $t/2 \leq |X| \leq |G|/2$, 
\item and any subgraph $F \subseteq G$ with $e(F) \leq d(G) \cdot \rho (|X|)|X|,$
\end{itemize}

\noindent then we have 
$$ |N_{G\backslash F}(X)| \geq \rho(|X|)|X|$$
\end{definition}

\noindent where $$\rho(x) = \rho(x,\epsilon_1,t) : =  \left\{ \begin{array}{lcl}
 0 & \mbox{if} & x < t/5, \\
\epsilon_1/\log^2(15x/t) & \mbox{if} & x \geq t/5.
\end{array}\right.$$\\

\noindent Note that the expansion rate of the expander above is only sublinear. Also, Haslegrave, Kim, and Liu~\cite{HKL20} gave the following statement that every graph contains one such sublinear  expander subgraph with almost the same average degree.\\

 \begin{theorem}[Haslegrave, Kim, and Liu~\cite{HKL20} ]\label{robust:thm} 
 There exists some $\epsilon_1> 0 $ such that the following holds for every $t>0$. Every 
   graph $G$  has an $(\epsilon_1,t)$-expander subgraph $H$ with $d(H) \geq d(G) / 2$ and 
  $\delta(H) \geq d(H) /  2$.\\
  \end{theorem}

\noindent A key property [Corollary 2.3~\cite{KS96}] of the  expanders  we use is the following lemma, that there exists a short path between two sufficiently large sets. This is formalized in the following statement.

\begin{lemma} [Haslegrave, Kim, and Liu~\cite{HKL20}] 
Let $\epsilon_1, t > 0$. If $G$ is an $n$-vertex $(\epsilon_1,t)$-expander, then any two vertex sets $X_1, X_2$, each of size at least $ x \leq t $, are of distance at most $m =\frac{1}{\epsilon_1} \log^3(15n/t)$ apart. This remains true even after deleting $\epsilon(x)\cdot x/4$ vertices from $G$.\\
\end{lemma}

\section{Proof of Theorem~\ref{maintheorem} }
In this section, we prove Theorem~\ref{maintheorem}. First, let us consider
 $$q = \sqrt{\frac{n}{d}}\ poly \log n {\text{ \ and \ }}k = \tilde{\Omega}\left(\frac{\sqrt{nd}}{\log^{10}n} \right),$$
\noindent where $k$ is the size of the complete minor we want to find.\\

We will construct $q+2$ disjoint sets $B_1, B_2, \cdots, B_q, D, U$ with $|B_i|=b= \frac{\epsilon^2 n}{k\log^{10}n}$ satisfying the following conditions:

\begin{itemize}
\item[(1)] $V(G) = D \cup B_1 \cup B_2 \cup \cdots B_q \cup U$
\item[(2)] for each $i \in [q]$, all induced subgraph $G[B_i]$ are connected and $|N(B_i)|\geq (1-2\epsilon)bd$
\item[(3)] $U$ is $(\epsilon/10, \alpha)$-expander
\item[(4)] $|N(D) \cap U | \ < \epsilon/10\cdot|D|$
\item[(5)] pairs $(B_i,B_j)$ satisfy that $B_i$ is connected to $B_j$ for at least $0.9q$ indices $j<i$
\item[(6)] $|D| \  \leq \frac{2}{\epsilon}|B|$, where $B = B_1 \cup B_2 \cup \cdots \cup B_q$.\\
\end{itemize}

\noindent Let us consider $q$ disjoint sets $ B_1, B_2, \cdots B_q $ with size  $|B_i|=b$ such that all induced subgraphs $G[B_i]$ are connected and $|N(B_i)|\ \geq (1-2\epsilon)bd$.\\

\noindent First, we claim that $D \leq 2 \epsilon n$. Since $N_G(D) \geq \epsilon |D|$, 
from  property $(4)$, we obtain
 $$|N_G(D) \cap B| = |N_G(D)| - |N_G(D)\cap U|\  \geq \ \epsilon |D| - \ \epsilon/10|D| \ = 9\epsilon/10|D|.$$
 
 \noindent Because of  $q\leq k$, we have $|B| \ \leq bq  \leq bk  \leq  \frac{\epsilon^2 n }{\log^{10} n} \leq \epsilon^2 n$.
 Therefore we  conclude that
 \begin{align}\label{equation:10}
 |D| \ \leq \ \frac{10}{9\epsilon}|B| \ \leq \frac{10}{9} \epsilon n\leq2\epsilon n.
 \end{align}
 
 \noindent It means that $|D|\ \leq \frac{2}{\epsilon}|B|$,  which is the property $(6)$.\\

 \noindent Let us consider the set $Bad =\{ \ i \leq q \ | \ | N_G(B_i) \cap U | \  < \frac{\epsilon}{10}\ bd\}$. \\
 
 \noindent Since $|N_G(B_i)|\ \geq (1-2\epsilon)\ bd$, for each $i \in Bad$, we derive that 
 \begin{align}\label{equation:11}
  e_G(B_i, D)  \geq |N_G(B_i)| - |N_G(B_i) \cap U| \ \geq (1-2\epsilon- \frac{\epsilon}{10})\ bd.\\ \nonumber
  \end{align}
 
 \noindent  
 To ensure the property $(5)$,  we suppose for the sake of contradiction that
 $|Bad| \ \geq 0.9q$.  
 Now we consider an average degree $d(G[\cup_{i\in Bad}B_i \cup D])$ on the induced  subgraph $G[\cup_{i\in Bad}B_i \cup D]$ as follows.\\
 
 \noindent Using the equations (\ref{equation:10}) and  (\ref{equation:11}), we have

 \begin{align*}
 d(G[\cup_{i\in Bad}B_i \cup D]) &\geq \frac{\sum_{q\in Bad} e_q (B_i, D)}{\sum_{i\in Bad}|B_i|+|D|}\\
 & \geq \frac{(1-3\epsilon)bd \cdot (0.9 q)}{bq+\frac{10}{9\epsilon}|B|} \\
 &\geq  \frac{(1-3\epsilon)0.9bdq}{\frac{2bq}{\epsilon}}\\
 &\geq  \frac{2\epsilon}{3}d 
 \end{align*}
\noindent which is a contradiction to the locally sparseness property.\\

\noindent For each $i \not \in Bad$, we consider $U_i= U \cap N_G(B_i)$ with the size $|U_i| \ \geq \frac{\epsilon}{10}\ bd \geq  \  \frac{\epsilon}{10}\sqrt{nd} \cdot  \log n$.
From the expansion of $G$, we obtain the following lemma that for a vertex in $X$ and integer $q = \sqrt{\frac{n}{d}}\cdot poly \log n$, there are at least  $\frac{\epsilon}{10}\sqrt{nd} \cdot  \log n$ vertices which are within distance at most $q = \sqrt{\frac{n}{d}}\cdot poly \log n$.\\

\begin{lemma}\label{lem:main}
 If we pick $\sqrt{nd} \cdot  \log n$ points $X$ at random, then we have
 $$ \sum_{i \not \in Bad} dist (X, U_i) \leq \sqrt{\frac{n}{d}}\cdot poly \log n.$$\\
\end{lemma}

 Using Lemma~\ref{lem:main}, we find a connected set $B'$ hitting all $N_G(B_i)$  with the size $|B_i| \  \leq \sqrt{\frac{n}{d}} \cdot \log^5 n$ for all $ i \not \in Bad$.
  Next, we  find a connected set $B''$ in $U$ with the size $|B''| = b - |B'|$  satisfying $N_G(B'') \geq \ (1-\epsilon/2)d \cdot |B''|$. If we let $B_{q+1}=B'\cup B''$, then we have $V(G) = D \cup B_1 \cup B_2 \cup \cdots B_{q+1} \cup U$ and $B = B_1 \cup B_2 \cup \cdots \cup B_q \cup B_{q+1}$.\\
  
  By the construction of the set $B_{q+1}$, we check that the property $(2)$ is satisfied. To ensure the property $(3)$, we move the following set $D_{q+1} \subseteq U$ to $D$. \\ 
  
  \newpage
  
  \noindent  Let us consider the maximal set $D_{q+1}$ such that 
  \begin{itemize}
  \item [(i)]$\frac{\epsilon}{10}$-expand in $U$
  \item [(ii)]$|D_{q+1}| \leq \epsilon n$\\
  \end{itemize}

   We suppose for the sake of contradiction that we can not move the set $D_{q+1}$ to $D$.  Let $U'$ be the largest $\frac{\epsilon}{2}$-expander in $U$.
   Let us consider the set $X$ which is not expander in $U'$.  Now we consider the set $D_{q+1} \cup X$. If $|X \cup D_{q+1}| \leq \epsilon n$, then the properties (i) and (ii) hold for the set $D_{q+1} \cup X $, which is a contradiction to the maximality of $D_{q+1}$. If $|X \cup D_{q+1}| > \epsilon n$, then we have
  
  $$ |N_G(X \cup D_{q+1})| \leq |B|+|D| + \frac{\epsilon n}{10} \leq \epsilon^4 n + 2 \epsilon^3 n + \frac{\epsilon n}{10}.$$\\

   \noindent This is a contradiction to the $\frac{\epsilon}{10}$-expand property of $U$. Therefore we move the set $D_{q+1}$ to $D$. Now we have the following $q+3$ disjoint sets $B_1, B_2, \cdots, B_q, B_{q+1}, D_{\text{new}}, U_{\text{new}}$ such that 
   
   \begin{itemize}
   \item[(a)] $B= B_1 \cup B_2 \cup \cdots, \cup B_{q} \cup B_{q+1}$
   \item[(b)] $ D_{\text{new}} = D_{\text{odd}} \cup D_{q+1}$
   \item[(c)]  $U_{\text{odd}} = U_{\text{new}} \cup D_{q+1} $.
   \end{itemize}
   
    From the expansion of $G$, we easily check that the property $(3)$ holds for the set $U_{\text{new}}$.  Since $|N(D_{q+1}) \cap U | \leq \frac{\epsilon}{10}|D_{q+1}|$, we observe that
    \begin{align*}
   \frac{\epsilon}{10}|D_{\text{odd}}| \ \geq |N(D_{\text{odd}}) \cap U_{\text{odd}}|     & \geq |N(D_{\text{new}}) \cap U_{\text{new}}| - |N(D_{q+1})\cap U_{\text{new}}| \\
    & \geq |N(D_{\text{new}}) \cap U_{\text{new}}| - \frac{\epsilon}{10}|D_{q+1}|.  
    \end{align*}
   
  \noindent  Therefore, we conclude that $$|N(D_{\text{new}}) \cap U_{\text{new}}| \leq \frac{\epsilon}{10}|D_{\text{new}}|,$$
  
  \noindent which is the property $(4)$.\\

    \noindent From property $(4)$, we obtain that 
     \begin{align*}
     |N_G(D_{\text{new}}) \cap B| & = |N_G(D_{\text{new}})| - |N_G(D_{\text{new}})\cap U_{\text{new}}| \\
     & \geq \epsilon \cdot |D_{\text{new}}| - \frac{\epsilon}{10}\cdot |D_{\text{new}}|  = \frac{9\epsilon}{10} \cdot |D_{\text{new}}|.
      \end{align*}
     
\noindent Then we conclude that 
\begin{align}\label{equation:n1}
|D_{\text{new}}| \leq \frac{10}{9\epsilon} |B| \leq \frac{2}{\epsilon} |B|
\end{align}
\noindent which is the  property $(6)$. \\
 
\noindent 
Let us consider the set $Bad =\{ \ i \leq q+1 \ | \ | N_G(B_i) \cap U_{\text{new}} | \  < \frac{\epsilon}{10}\ bd\}$. \\

 \noindent By the construction of $B_{q+1}$,   for each $i \in Bad$ we have $|N_G(B_i)|\ \geq (1-2\epsilon)\ bd$. \noindent From property $(4)$, we derive that \\
 \begin{align}\label{equation:n2}
  e_G(B_i, D_{\text{new}})  \geq |N_G(B_i)| - |N_G(B_i) \cap U_{\text{new}}| \ \geq (1-2\epsilon- \frac{\epsilon}{10})\ bd.\\\nonumber
  \end{align}
 
 \noindent To ensure  the property $(5)$, we suppose for the sake of contradiction that
 $|Bad| \ \geq 0.9q$.  
 Now we consider an average degree $d(G[\cup_{i\in Bad}B_i \cup D_{\text{new}}])$ on the induced  subgraph $G[\cup_{i\in Bad}B_i \cup D_{\text{new}}]$ as follows.\\
 
 \noindent Using the equations (\ref{equation:n1}) and  (\ref{equation:n2}), we have

 \begin{align*}
 d(G[\cup_{i\in Bad}B_i \cup D_{\text{new}}]) &\geq \frac{\sum_{q\in Bad} e_q (B_i, D_{\text{new}})}{\sum_{i\in Bad}|B_i|+|D_{\text{new}}|} \geq \frac{(1-3\epsilon)bd \cdot (0.9 q)}{bq+\frac{10}{9\epsilon}|B|} \\
 &\geq  \frac{(1-3\epsilon)0.9bdq}{\frac{2bq}{\epsilon}}\geq  \frac{2\epsilon}{3}d \\
 \end{align*}
\noindent which is a contradiction to the locally sparseness property. We repeat this process as long as $bq \leq \frac{\epsilon^{2} n}{\log^{10} n}.$\\

{\bf{Acknowledgement.}} I would like to thank Jaehoon Kim and Hong Liu for helpful discussions.

\end{document}